\documentclass{amsart}

\usepackage{amsmath}
\usepackage{amsthm}
\usepackage{amsfonts}
\usepackage{amssymb}

\usepackage{graphics}
\usepackage{epsfig}
\usepackage{color}
\usepackage{verbatim}

\newcommand{\R}{{\mathbb R}}
\newcommand\BL{\operatorname{BL}}
\newcommand{\datum}{\mathbf}
\newcommand{\codim}{\operatorname{codim}}
\newcommand{\maybed}{}
\newcommand{\maybedsigma}{S}

\theoremstyle{plain}
  \newtheorem{thm}{Theorem}[section]

  \newtheorem{prp}[thm]{Proposition}

\theoremstyle{definition}
  
  \newtheorem*{exam}{Example}
  \newtheorem*{ack}{Acknowledgement}

\numberwithin{equation}{section}

\title{Higher order transversality in harmonic analysis}
\begin{document}
\author{Jonathan Bennett}
\address{School of Mathematics, The Watson Building, University of Birmingham, Edgbaston,
Birmingham, B15 2TT, England}
\email{j.bennett@bham.ac.uk}
\thanks{Supported by JSPS Kakenhi no. 19H01796 (Bez)}

\author{Neal Bez}
\address{Department of Mathematics, Graduate School of Science and Engineering,
Saitama University, Saitama 338-8570, Japan}
\email{nealbez@mail.saitama-u.ac.jp}



\begin{abstract} 
In differential topology two smooth submanifolds $S_1$ and $S_2$ of euclidean space are said to be transverse if the tangent spaces at each common point together form a spanning set.
The purpose of this article is to explore a much more general notion of transversality pertaining to a \emph{collection} of submanifolds of euclidean space. In particular, we show that three seemingly different concepts of transversality arising naturally in harmonic analysis, are in fact equivalent. This result is an amalgamation of several recent works on variants of the Brascamp--Lieb inequality, and we take the opportunity here to briefly survey this growing area. This is not intended to be an exhaustive account, and the choices made reflect the particular perspectives of the authors.
\end{abstract}
\maketitle

\section{Introduction}
In differential topology, two smooth submanifolds $S_1, S_2$ of $\mathbb{R}^n$ are said to be \emph{transverse} if at each point of their intersection, the tangent spaces of $S_1, S_2$ together span $\mathbb{R}^n$ -- that is
\begin{equation}\label{trandiff}
x\in S_1\cap S_2\implies T_xS_1+T_xS_2=\mathbb{R}^n.
\end{equation}
Of course $(T_xS_1+T_xS_2)^\perp=(T_xS_1)^\perp\cap (T_xS_2)^\perp$, and so this notion of transversality is equivalent to requiring that the normal spaces to $S_1, S_2$ at a common point intersect trivially.

One reason why transversality is important in analysis is that it allows us to make sense of products of distributions. For example, if $S_1, S_2$ are transverse smooth submanifolds of $\mathbb{R}^n$, and for each $j=1,2$ the $S_j$-carried distribution $\delta_{S_j}$ is given by
$$
\langle\delta_{S_j},\varphi\rangle=\int_{S_j}\varphi\:d\sigma_j,
$$
where $d\sigma_j$ denotes surface measure on $S_j$, then the product distribution $\delta_{S_1}\delta_{S_2}$ is well-defined and given by
$$
\langle\delta_{S_1}\delta_{S_2},\varphi\rangle=\int_{S_1\cap S_2}\varphi \:d\mu,
$$
for a measure $d\mu$ that is absolutely continuous with respect to surface measure on $S_1\cap S_2$ (see, for example, Sogge \cite{Sogge} and Foschi--Oliveira e Silva \cite{FosOli}).

Assigning meaning to such distributional products becomes relevant in harmonic analysis when defining convolutions of measures carried on submanifolds of $\mathbb{R}^n$ -- something that arises frequently. In this context it is natural to strengthen \eqref{trandiff} in a way that makes it translation-invariant in $S_1$ and $S_2$ independently, so that it becomes \begin{equation}\label{convtrans}
x_1\in S_1, x_2\in S_2\implies T_{x_1}S_1+T_{x_2}S_{2}=\mathbb{R}^n.
\end{equation}
With this translation-invariant notion of transversality, it follows that the convolution $d\sigma_1*d\sigma_2$ is absolutely continuous with respect to Lebesgue measure on $\mathbb{R}^n$. Moreover, if $S_1, S_2$ are compact, then $d\sigma_1*d\sigma_2$ has bounded density, or in other words, the (bilinear) estimate
\begin{equation}\label{conv2}
\|g_1d\sigma_1*g_2d\sigma_2\|_{L^\infty(\mathbb{R}^n)}\lesssim\|g_1\|_{L^\infty(S_1)}\|g_2\|_{L^\infty(S_2)}
\end{equation}
holds.

Such convolutions of surface carried measures arise naturally in harmonic analysis, particularly in the restriction theory of the Fourier transform. This theory concerns the $L^q(\mathbb{R}^n)$ integrability properties of Fourier transforms of measures with $L^p$ densities supported on (typically curved) submanifolds of $\mathbb{R}^n$, and has many applications, from dispersive partial differential equations to analytic number theory (see, for example, Stovall \cite{StovRest}). Questions in this area are usually formulated in terms of the \emph{extension operator} $g\mapsto\widehat{gd\sigma}$, given by
$$
\widehat{gd\sigma}(x)=\int_{S}e^{i\langle x , \xi \rangle}g(\xi)d\sigma(\xi),
$$
where $S$ is a smooth (typically compact) submanifold of $\mathbb{R}^n$, and $d\sigma$ is surface measure on $S$, as before. At the centre of this theory is the celebrated \emph{restriction conjecture} of Stein, which states that if $S$ is a smooth compact hypersurface with everywhere nonvanishing gaussian curvature, then
$$
\|\widehat{gd\sigma}\|_{L^q(\mathbb{R}^n)}\lesssim\|g\|_{L^p(S)}
$$
provided $\tfrac{1}{q}<\tfrac{n-1}{2n}$ and $\tfrac{1}{q}\leq\tfrac{n-1}{n+1}\tfrac{1}{p'}$ (see Stein \cite{bigStein}). This conjecture is settled for $n=2$  (C. Fefferman and Stein \cite{Feff}, \cite{bigStein}; see also Zygmund \cite{Zyg}), and is still open for $n\geq 3$ -- see Hickman--Rogers \cite{HR} for further discussion and recent developments.
The relevance of transversality in restriction theory stems from the simple observation that if $S_1, S_2$ are suitable pieces of such a curved manifold $S$, then they will be transverse in the sense of \eqref{convtrans}. For example, this will be the case if $S$ is a hemisphere and $S_1, S_2\subseteq S$ are disjoint caps. This naturally leads one to consider bilinear extension operators (or bilinear interactions) of the form
\begin{equation}\label{convprod}
(g_1,g_2)\mapsto \widehat{g_1d\sigma_1}\widehat{g_2d\sigma_2}=(g_1d\sigma_1*g_2d\sigma_2)\:\widehat{\:},
\end{equation}
and seek estimates of the form
\begin{equation}\label{bil2}
\|\widehat{g_1d\sigma_1}\widehat{g_2d\sigma_2}\|_{L^{q/2}(\mathbb{R}^n)}\lesssim\|g_1\|_{L^{p_1}(S_1)}\|g_2\|_{L^{p_2}(S_2)},
\end{equation}
under the assumption that $S_1, S_2$ are transverse. We refer to Bourgain \cite{Bourgain93} and Tao--Vargas--Vega \cite{TVV} for the origins of this idea.

Under a transversality hypothesis alone, that is, not stipulating any curvature properties of $S_1, S_2$, the bilinear estimate \eqref{bil2} is well understood, and the whole story may be reduced to the endpoint estimate
\begin{equation}\label{rest2}
\|\widehat{g_1d\sigma_1}\widehat{g_2d\sigma_2}\|_{L^2(\mathbb{R}^n)}\lesssim\|g_1\|_{L^2(S_1)}\|g_2\|_{L^2(S_2)}.
\end{equation}

In two dimensions ($n=2$) matters are particularly simple, and the statements \eqref{convtrans}, \eqref{conv2} and \eqref{rest2} are easily seen to be \emph{equivalent}, in that 
$$
\eqref{convtrans} \implies \eqref{conv2} \implies\eqref{rest2}\implies\eqref{convtrans}
$$
may be quickly verified. The first implication has already been discussed, and the second amounts to an application of Plancherel's theorem via \eqref{convprod}, followed by a routine interpolation argument. The third follows from the standard (Knapp-type) examples in this context -- see Section \ref{SectProof}.

The main purpose of this article is to present a broad generalisation of the simple equivalence above. Our main result is the following, which unifies the recent works Bennett--Carbery--Christ--Tao \cite{BCCT1, BCCT2}, Bennett--Bez--Flock--Lee \cite{BBFL} and Bennett--Bez--Buschenhenke--Cowling--Flock \cite{BBBCF}.
\begin{thm}\label{tfae}
Suppose $S_1,\hdots,S_m$ are smooth compact submanifolds of $\mathbb{R}^n$, and $p_1,\hdots, p_m \in [1,\infty]$ are Lebesgue exponents satisfying
\begin{equation}\label{scaleawayscaleawayscaleaway}
\sum_{j=1}^m \frac{\dim(S_j)}{p_j'}=n.
\end{equation}
Then the following are equivalent:
\begin{equation}\tag{C}
\|g_1d\sigma_1*\cdots *g_md\sigma_m\|_{L^\infty(\mathbb{R}^n)}\lesssim \|g_1\|_{L^{p_1}(S_1)}\cdots\|g_m\|_{L^{p_m}(S_m)}
\end{equation}
\begin{equation}\tag{E}
\int_{B_R}\prod_{j=1}^m|\widehat{g_jd\sigma_j}|^{2/p_j'}\lesssim_\varepsilon R^\varepsilon\prod_{j=1}^m\|g_j\|_{L^2(S_j)}^{2/p_j'}\;\;\;\mbox{ all }\varepsilon>0
\end{equation}
\begin{equation}\tag{T}%\tag{``$\mathbf{p}$-transversality"}\label{gent}
\dim(V)\leq \sum_{j=1}^m\frac{\dim(V)-\dim(V\cap(TS_j)^\perp)}{p_j'}
\end{equation}
for all subspaces $V$ of $\mathbb{R}^n$ and all tangent spaces $TS_j$ of $S_j$.
\end{thm}
The implications (T)$\implies$(E)  and (T)$\implies$(C) were proved at this level of generality in \cite{BBFL} and \cite{BBBCF} respectively, both using the method of induction-on-scales. That (C)$\implies$(T) follows from Bennett--Bez--Guti\'errez \cite{BBG}, and that (E)$\implies$(T) will be presented in Section \ref{SectProof}.
Some words on nomenclature: here (C), (E) and (T) stand for ``convolution", ``extension" and ``transversality" respectively. Of course this notion of transversality depends on the exponents $\mathbf{p}=(p_1,\hdots, p_m)$, and it should be noticed that if $m=2$ and $p_1=p_2=\infty$, then it coincides with the familiar \eqref{convtrans}. 

\begin{exam} 
Suppose $m=n$, the submanifolds $S_1,\hdots,S_n$ are \emph{hypersurfaces}, and $p_1=\cdots = p_n= (n-1)'$. In this case the transversality condition (T) amounts to the statement that 
\begin{equation}\label{lwt} (TS_1)^\perp+\cdots +(TS_n)^\perp=\mathbb{R}^n;\end{equation}
in other words, any selection of unit normals $\nu_1,\hdots,\nu_n$ to $S_1,\hdots,S_n$ respectively, forms a basis of $\mathbb{R}^n$
(note that the volume form $|\nu_1\wedge\cdots\wedge \nu_n|$ is also bounded below by compactness). In this case the extension estimate (E) becomes the well-known multilinear restriction inequality 
\begin{equation}\label{mlkbct}
\|\widehat{g_1d\sigma_1}\cdots\widehat{g_nd\sigma_n}\|_{L^{\frac{2}{n-1}}(B_R)}\lesssim_\varepsilon R^\varepsilon\|g_1\|_{L^2(S_1)}\cdots\|g_n\|_{L^2(S_n)}
\end{equation}
of Bennett--Carbery--Tao \cite{BCT}. When $n=3$ the estimate (C) amounts to the statement that the convolution of $L^2$ densities supported on $S_1$ and $S_2$ restricts to an $L^2$ density on $S_3$ -- see Bejenaru--Herr--Tataru \cite{BHT}, and Bennett--Bez \cite{BB} for (C) in higher dimensions. Underlying this is a certain nonlinear perturbation of the classical Loomis--Whitney inequality of Bennett--Carbery--Wright \cite{BCW}. In the context of the particular transversality condition \eqref{lwt} the inequalities (E) and (C) have had numerous applications -- see Section \ref{Sect:App} for some examples.
\end{exam}

Some remarks on the condition \eqref{scaleawayscaleawayscaleaway} are in order. This condition ensures that the estimates (E) and (C) are ``curvature blind", and it has this effect by ensuring that they are both\footnote{Strictly speaking this scale-invariance requires that $\varepsilon=0$ in (E).} invariant under isotropic scalings of the submanifolds $S_j$ (with respect to scalings of the underlying euclidean space). Of course curvatures may be made arbitrarily small under isotropic dilations, whereas quantitative measures of transversality are left unchanged. As an example, the well-known three-dimensional bilinear extension estimate 
\begin{equation}\label{tt}
\|\widehat{g_1d\sigma_1}\widehat{g_2d\sigma_2}\|_{L^{5/3}(B_R)}\lesssim_\varepsilon R^\varepsilon\|g_1\|_{L^2(S_1)}\|g_2\|_{L^2(S_2)}
\end{equation}
of Tao \cite{TaoBil}, involving separated compact subsets $S_1,S_2$ of the paraboloid in $\mathbb{R}^3$ fails to satisfy the scaling condition \eqref{scaleawayscaleawayscaleaway}. However, despite \eqref{tt} being (E) with $p_1=p_2=6$, the transversality condition (T) with the same exponents is easily seen to fail. The point here is that the curvature of the paraboloid is playing a role in \eqref{tt}, along with the transversality. We note the consistency of this observation with our assertions about \eqref{rest2}.

It is natural to conjecture that (E) holds with $\varepsilon=0$, as is done in \cite{BCT} for the special case \eqref{mlkbct}. The power loss may at least be reduced to polylogarithmic loss, replacing the $R^\varepsilon$ factor with a power of $\log R$, as is observed in Bennett \cite{B} and Zhang \cite{Z} for \eqref{mlkbct} and (E) respectively, although removing it entirely is only currently possible in degenerate or very simple cases, such as when $S_j=\mathbb{R}^n$ for each $j$, or when $n=2$. Very recently however, this loss has been successfully removed from all proper interpolants of \eqref{mlkbct} and the elementary 
bound $$
\|\widehat{g_1d\sigma_1}\cdots\widehat{g_nd\sigma_n}\|_{L^\infty(\mathbb{R}^n)}\lesssim\|g_1\|_{L^1(S_1)}\cdots\|g_n\|_{L^1(S_n)}$$
%$L^1(S_1)\times\cdots\times L^1(S_n)\rightarrow L^\infty(\mathbb{R}^n)$ bound 
-- see Tao \cite{Tao19}. Prior to \cite{Tao19} such global estimates away from the sharp line had been achieved under the additional hypothesis that the hypersurfaces $S_j$ have everywhere nonvanishing gaussian curvature, using $\varepsilon$-removal techniques -- see Bourgain--Guth \cite{BG}.

The key to understanding Theorem \ref{tfae} lies in the theory of the Brascamp--Lieb inequality, first formulated by Brascamp and Lieb \cite{BL} as a generalised form of Young's convolution inequality. As we shall see, the three statements (C), (E) and (T) are manifestations of certain seemingly quite different generalisations of this inequality, whose equivalence has only recently been understood.\footnote{Notwithstanding the expectation that (E) holds with $\varepsilon=0$.}

\section{The Brascamp--Lieb inequality and its many variants and generalisations}\label{SectWider}
The Brascamp--Lieb inequality is a well-known and far-reaching generalisation of a wide range of sharp functional inequalities in analysis, including the multilinear H\"older, Loomis--Whitney and Young convolution inequalities.
It takes the form
\begin{equation}\label{BL}
\Biggl|\int_{\R^n}\prod_{j=1}^m f_j(L_jx) \,dx\Biggr|
\leq \BL(\datum{L},\datum{p})\prod_{j=1}^m\|f_j\|_{L^{p_j}(\R^{n_j})},
\end{equation}
where the mappings $L_j:\R^n\to\R^{n_j}$ are linear surjections, $p_j\in[1,\infty]$, and $\BL(\datum{L},\datum{p})$ denotes the smallest constant (which may be infinite).
We refer to $$(\datum{L},\datum{p})=((L_j)_{j=1}^m,(p_j)_{j=1}^m)$$ as the \textit{Brascamp--Lieb datum}, and $\BL(\datum{L},\datum{p})$ as the \textit{Brascamp--Lieb constant}.\footnote{While the functions $f_j$ may be complex valued here, in studying the Brascamp--Lieb constant $\BL(\mathbf{L},\mathbf{p})$ we may of course restrict attention to nonnegative $f_j$, and this is often implicit. We caution that the Brascamp--Lieb datum is usually presented in terms of the reciprocals of the exponents $p_j$, as elements of $(0,1]$ -- see \eqref{BLL}.}
This inequality was first  formulated in \cite{BL} and gives rise to an elegant and powerful theory; further notable contributions in this direction include those of Ball \cite{Ball}, Lieb \cite{L}, Barthe \cite{Barthe}, Carlen, Lieb and Loss \cite{CLL}, as well as \cite{BCCT1} and Barthe, Cordero-Erausquin, Ledoux and Maurey \cite{BCLM}.
Ball \cite{Ball} derived a particular class of rank-one Brascamp--Lieb inequalities and pioneered their use in convex geometry, solving several problems on volumes of sections of convex bodies.
Barthe \cite{Barthe, BartheICM} extended such geometric Brascamp--Lieb inequalities to the general rank case, and introduced the use of optimal transport methods in order to advance the general theory of the Brascamp--Lieb inequality and its dual form.
Applications and perspectives on the Brascamp--Lieb inequality may be found even more widely, including in convex geometry \cite{BGMN, Braz, BrazG, Gardner}, probability, stochastic processes and statistics \cite{ALS, ALS2, Finner, Lehec}, information theory \cite{BCM, CC, CLL, EFKY, LCCV}, scattering theory \cite{AFR, Brown, BOP, Perry}, combinatorics \cite{BT, ET}, group theory \cite{CLL2, Eberhard}, and theoretical computer science \cite{CDKSY0, CDKSY, DGOS, GGOW, VY}. For example, the information-theoretic connection reveals that the Brascamp--Lieb inequality has an equivalent formulation in terms of a generalised notion of subadditivity of the entropy (see \cite{CC}).

Before introducing some of the recent generalisations of \eqref{BL}, we briefly describe some of the key features and examples of the classical theory. This is not intended to be an exhaustive or balanced account, and the reader is referred to \cite{BCCT1, BCCT2} for further results and discussion.

A well-known example, and indeed the historical motivation for the Brascamp--Lieb inequality, is the celebrated sharp version of Young's convolution inequality, first proved by Beckner \cite{Beckner2, Beckner} and Brascamp--Lieb \cite{BL}.
In the framework of \eqref{BL}, this may be stated as
\begin{equation}\label{Young}
\begin{aligned}
\int_{\R^d} \int_{\R^d}f_1(y)f_2(x-y)f_3(x) \,dy \,dx
\leq (C_{p_1}C_{p_2}C_{p_3})^d
\|f_1\|_{L^{p_1}(\mathbb{R}^d)}\|f_2\|_{L^{p_2}(\mathbb{R}^d)}\|f_3\|_{L^{p_3}(\mathbb{R}^d)},
\end{aligned}
\end{equation}
where $p_1,p_2,p_3\in [1,\infty]$, $\tfrac{1}{p_1}+\tfrac{1}{p_2}+\tfrac{1}{p_3}=2$, and $C_r=((1-1/r)^{1-1/r}/(1/r)^{1/r})^{1/2}$.
The main significance of this statement lies in the conclusion that, on the relative interior of the set of admissible exponents $(\tfrac{1}{p_1},\tfrac{1}{p_2},\tfrac{1}{p_3})$, the optimal constant $(C_{p_1}C_{p_2}C_{p_3})^d$ is strictly less than $1$, and is uniquely attained on suitably scaled isotropic centred gaussian inputs $f_j$.   Remarkably, this phenomenon turns out to be quite typical in the general context of \eqref{BL}. This is captured by a fundamental theorem of Lieb \cite{L}, which guarantees the existence of extremising sequences of centred gaussians for \eqref{BL}. This reduces the complexity of working with the Brascamp--Lieb constant considerably since it quickly gives rise to the formula
$$
\BL(\mathbf{L},\mathbf{p})=\sup_{\mathbf{A}}\frac{\prod_{j=1}^m\det(A_j)^{1/(2p_j)}}{\det\Bigl(\sum_{j=1}^m \frac{1}{p_j}L_j^*A_jL_j\Bigr)^{1/2}},
$$
where the supremum is taken over all $m$-tuples $\mathbf{A}=(A_1,\hdots, A_m)$ of positive definite $n_j\times n_j$ symmetric matrices $A_j$, $1\leq j\leq m$. In particular, the proof of (T)$\implies$(C) in Theorem \ref{tfae} relies heavily on a quantified version of this fact -- see \cite{BBBCF}. We refer the reader to \cite{BCCT1} and the references there for further structural results, such as statements on the existence and uniqueness of gaussian extremisers.

From the point of view of the restriction theory of the Fourier transform, arguably the most important examples are of ``Loomis--Whitney" type (we refer the reader forward to Section \ref{Sect:ext} for elaboration of this). Here the kernels of the linear surjections $L_j:\mathbb{R}^n\rightarrow\mathbb{R}^{n_j}$ collectively form a basis for $\mathbb{R}^n$ -- that is
\begin{equation}\label{LoomW}
\ker(L_1)\oplus\cdots\oplus\ker(L_m)=\mathbb{R}^n.
\end{equation}
For such data simple examples reveal that $\BL(\mathbf{L},\mathbf{p})<\infty$ if and only if $p_1=\cdots=p_n=m-1$, and in this case there is an explicit expression for
$\BL(\mathbf{L},\mathbf{p})$ in terms of the natural volume form associated with $\mathbf{L}$ -- see for example \cite{BB}.

The inequality \eqref{BL} may be interpreted as a bound on multilinear forms of the type
\begin{equation}\label{stefanlikes}
(f_1,\dots,f_m)\mapsto\int_H f_1\otimes\dots\otimes f_m \,d\mu_H,
\end{equation}
where $H$ is a subspace of the cartesian product $\mathbb{R}^{n_1}\times\cdots\times\mathbb{R}^{n_m}$, and integration is with respect to Lebesgue measure on $H$.
Here $H$ encodes the linear maps $\datum{L}$ as the range of $x \mapsto (L_1x,\dots,L_mx)$, allowing the pair $(H,\mathbf{p})$ to be interpreted as a certain parametrisation-free Brascamp--Lieb datum, and $\BL(H,\mathbf{p})$, defined to be the best constant in the inequality
\begin{equation}\label{BLPF}
\Biggl|\int_H f_1\otimes\dots\otimes f_m \,d\mu_H\Biggr|\leq \BL(H,\datum{p})\prod_{j=1}^m\|f_j\|_{L^{p_j}(\R^{n_j})},
\end{equation}
a parametrisation-free Brascamp--Lieb constant.\footnote{While this is clearly an abuse of notation, the context here eliminates any possible confusion.}
As one may expect, quite how this subspace $H$ sits relative to the coordinate subspaces $\mathbb{R}^{n_1},\hdots,\mathbb{R}^{n_m}$ determines the finiteness (or otherwise) of the Brascamp--Lieb constant, and it is this that is ultimately captured by the transversality condition (T) in the statement of Theorem \ref{tfae}. The following finiteness characterisation will be crucial in making this connection.
\begin{thm}[\cite{BCCT1}]\label{t:fin}
The Brascamp--Lieb constant $\BL(\mathbf{L},\mathbf{p})$ is finite if and only if
\begin{equation}\label{sca}
n=\sum_{j=1}^m\frac{n_j}{p_j}
\end{equation}
and
\begin{equation}\label{scaV}
\dim(V)\leq\sum_{j=1}^m\frac{\dim(L_jV)}{p_j}
\end{equation}
for all subspaces $V$ of $\mathbb{R}^n$.
\end{thm}
We remark that Theorem \ref{t:fin} follows similar results in \cite{Barthe}, and the case where the maps $L_j$ have rank one was established in \cite{CLL}.

In terms of the parametrisation-free data $(H,\mathbf{p})$, the above result may be re-interpreted as follows: $\BL(H,\mathbf{p})$ is finite if and only if
\begin{equation}\label{scapf}
\dim(H)=\sum_{j=1}^m\frac{n_j}{p_j}
\end{equation}
and
\begin{equation}\label{scapfV}
\dim(V)\leq\sum_{j=1}^m\frac{\dim(\pi_jV)}{p_j}
\end{equation}
for all subspaces $V$ of $H$. Here $\pi_j$ denotes the orthogonal projection from $\mathbb{R}^{n_1}\times\cdots\times\mathbb{R}^{n_m}$ onto the $j$th factor $\mathbb{R}^{n_j}$.

There is one clear advantage to the parametrisation-free formulation \eqref{BLPF} stemming from the elementary fact that the structure of the multilinear form involved is manifestly \emph{Fourier--invariant} -- that is,
$$
\int_H f_1\otimes\dots\otimes f_m \,d\mu_H=\int_{H^\perp} \widehat{f}_1\otimes\dots\otimes \widehat{f}_m \,d\mu_{H^\perp},\footnote{Of course nonnegativity is not Fourier-invariant, so it is important that we consider general complex-valued functions $f_j$ in this context.}
$$
where $H^\perp$ denotes the orthogonal complement of $H$ in $\mathbb{R}^{n_1}\times\cdots\times\mathbb{R}^{n_m}$. This gives rise to a useful (Fourier-) duality principle, which states that
\begin{equation}\label{inv}
\BL(H,\mathbf{p})=C_{\mathbf{p},\mathbf{n}}\BL(H^\perp,\mathbf{p}');
\end{equation}
see \cite{BBBCF} (and also Bennett--Jeong \cite{BJ} for a discrete analogue). Here $\mathbf{p}'=(p_1',\hdots,p_m')$, and $C_{\mathbf{p},\mathbf{n}}$ denotes an explicit constant depending on $\mathbf{p}$ and the underlying dimensions $\mathbf{n}=(n_1,\hdots,n_m)$. This plays an important role in establishing that (T)$\implies$(C) in the proof of Theorem \ref{tfae} -- see Section \ref{SectProof}.

Finally, we note that the Brascamp--Lieb inequality is usually stated in the form
\begin{equation}\label{BLL}
\int_{\mathbb{R}^n}\prod_{j=1}^mf_j(L_jx)^{r_j}dx\leq \widetilde{\BL}(\mathbf{L},\mathbf{r})\prod_{j=1}^m\Bigl(\int_{\mathbb{R}^{n_j}}f_j\Bigr)^{r_j},
\end{equation}
where $r_j\in [0,1]$, and the functions $f_j$ are nonnegative and integrable on their respective domains. On replacing $f_j$ in \eqref{BL} by $f_j^{r_j}$, where $r_j=1/p_j$, we see that \eqref{BL} (with nonnegative $f_j$) and \eqref{BLL} are the same inequality, and in particular $\widetilde{\BL}(\mathbf{r},\mathbf{L})=\BL(\mathbf{p},\mathbf{L})$. The advantage of this formulation is that it alludes to a certain self-similarity structure that may be exploited by the method of induction-on-scales -- see, for example, \cite{BB} or \cite{BBBCF} for further discussion. In order to avoid confusion we make no further reference to $\widetilde{\BL}$.

In recent years a number of variants of the Brascamp--Lieb inequality have emerged in harmonic analysis, with a range of applications. This effectiveness reflects an improved understanding of the role of curvature in harmonic analysis, and in particular its relation to transversality of the type discussed in this article. In the remainder of this section we survey a variety of interconnected set-ups, with the first and third of particular importance for the proof of Theorem \ref{tfae}.

\subsection{A nonlinear variant}
The Brascamp--Lieb inequality, stated in either form \eqref{BL} or \eqref{BLPF}, involves linearity in an important way -- either the linearity of the maps $L_j$, or the linearity of the subspace $H$. A variety of problems in harmonic analysis and PDE raise questions about the necessity of this linear structure. For example, may the $L_j$ be replaced by smooth submersions $B_j$, or equivalently, may $H$ be replaced by a smooth submanifold $M$, at least in a neighbourhood of a point of $\mathbb{R}^n$? Affirmative answers to these questions are known as \emph{nonlinear Brascamp--Lieb inequalities}, and this line of research traces back to \cite{BCW}. At the level of local statements, the following near-optimal result was obtained recently in \cite{BBBCF}:
\begin{thm}[\cite{BBBCF}]\label{main'}
Suppose $M$ is a $C^2$ submanifold of $\mathbb{R}^{n_1}\times\cdots\times\mathbb{R}^{n_m}$ and $x\in M$. Then given any $\varepsilon>0$ there exists $\delta>0$ such that
\begin{equation}\label{great!}
\int_{M\cap B(x,\delta)}f_1\otimes\cdots\otimes f_m \,d\mu_M\leq (1+\varepsilon)\BL(T_xM,\mathbf{p})\prod_{j=1}^m\|f_j\|_{L^{p_j}(\R^{n_j})}.
\end{equation}
Here $d\mu_M$ denotes Lebesgue measure on $M$, and $T_xM$ denotes the tangent space to $M$ at $x$.
\end{thm}
Theorem \ref{main'} follows several nonlinear Brascamp--Lieb inequalities established under various additional (structural) hypotheses on the manifold $M$, and obtained by different methods -- see \cite{BCW, BHT, BB, KS, CHV, BBBCF} for further discussion. We also refer the reader to recent developments by Duncan \cite{Duncan1, Duncan2} on nonlinear Brascamp--Lieb inequalities, including certain global estimates and  stability results.

It is perhaps worth remarking that many positive multilinear forms way be expressed as 
\begin{equation}\label{functionalanalysis}
(f_1,\dots, f_m)\mapsto\int_{\R^{n_1}\times\dots\times\R^{n_m}} f_1\otimes\dots\otimes f_m \,d\mu
\end{equation}
for an appropriate measure $\mu$. Indeed, the Schwartz kernel theorem expresses any $m$-linear form on $m$-tuples of Schwartz functions $(f_1,\hdots,f_m)\in\mathcal{S}(\mathbb{R}^{n_1})\times\cdots\times\mathcal{S}(\mathbb{R}^{n_m})$, acting continuously in each component, as an element of  $\mathcal{S}'(\mathbb{R}^{n_1+\cdots+n_m})$ acting on the tensor product $f_1\otimes\cdots\otimes f_m$. If this distribution extends to a continuous linear functional on $C_0$, then 
a version of the Riesz representation theorem yields the representation \eqref{functionalanalysis} for some Radon measure $\mu$. Theorem \ref{main'} provides Lebesgue space bounds on \eqref{functionalanalysis} in situations where $\mu$ specialises to integration over a submanifold $M$ satisfying certain structural (transversality) and regularity hypotheses.
Multilinear forms of this type, which are sometimes referred to as \textit{multilinear Radon-like transforms}, arise frequently in analysis -- see Section \ref{Sect:App} for some examples in the context of dispersive PDE.

\subsection{A Kakeya-type variant}\label{Sect:kak}
Informally speaking Kakeya-type problems concern the extent to which families of geometric objects (usually subsets of $\mathbb{R}^n$) may be arranged so as to minimise the space in which they occupy (or maximise the extent to which they overlap). Usually the geometric objects are $\delta$-neighbourhoods of unit line segments $T\subseteq\mathbb{R}^n$, referred to as $\delta$-tubes, belonging to some family $\mathbb{T}$, and the objective is to control the ``multiplicity function"
$$
x\mapsto\sum_{T\in\mathbb{T}}\chi_T(x),
$$
where $\chi_T$ denotes the characteristic function of $T$. For families $\mathbb{T}$ of tubes whose direction set forms a $\delta$-separated subset of the unit sphere $\mathbb{S}^{n-1}$, the Kakeya maximal conjecture states that for every $\varepsilon>0$,
\begin{equation}\label{lkak}
\Bigl\|\sum_{T\in\mathbb{T}}\chi_T\Bigr\|_{L^{\frac{n}{n-1}}(\mathbb{R}^{n})}\lesssim_{\varepsilon}\delta^{-\varepsilon}(\delta^{n-1}\#\mathbb{T})^{\frac{n-1}{n}}
\end{equation}
uniformly in the family $\mathbb{T}$. A notable consequence of this is the Kakeya set conjecture, which asserts that a Kakeya set (a compact subset of $\mathbb{R}^n$ containing a unit line segment in every direction) has full Hausdorff dimension -- a statement that is only known in two dimensions; see for example Hickman--Rogers--Zhang \cite{HRZ} or Zahl \cite{Zahl} for some historical background and the current status of this active problem in higher dimensions.

The Kakeya maximal conjecture is well known to follow from the restriction conjecture (stated in the introduction) -- see \cite{BCSS}. The mechanism by which this follows relates the (unit) normal vectors to the hypersurface $S$ with the directions within the family of tubes, and through this the $\delta$-separation condition may be viewed as a sort of curvature condition. As we shall see, it is also natural to consider ``multilinear" variants of the Kakeya conjecture \eqref{lkak} where the $\delta$-separation condition is replaced with a suitable transversality condition between \emph{families} of tubes.
A very general Kakeya-type inequality of this flavour, referred to as the \emph{Kakeya--Brascamp--Lieb inequality}, was first established in Zhang \cite{Z} (see also Zorin--Kranich \cite{ZK}), following a weaker version in \cite{BBFL} involving an $\varepsilon$-loss similar to that in \eqref{lkak}. Both of these were preceded by special cases corresponding to Loomis--Whitney-type data, established by Guth \cite{Guth} and Bennett--Carbery--Tao \cite{BCT} respectively -- see also Guth \cite{Guth2}. We also refer to Bourgain--Guth \cite{BG} and Carbery--Valdimarsson \cite{CV} for further results.

In order to state a form of the Kakeya--Brascamp--Lieb inequality let $V_1,\hdots, V_m$ denote subspaces on $\mathbb{R}^n$ of codimensions $n_1,\hdots, n_m$ respectively, and suppose that for each $1\leq j\leq m$ the set $\mathbb{T}_j$ denotes a finite collection of $\delta$-neighbourhoods of codimension-$n_j$ affine subspaces that are, modulo translations, sufficiently close to $V_j$ (with respect to the standard metric on the grassmann manifold of codimension-$n_j$ subspaces of $\mathbb{R}^n$). The following is a special case of a more general result of Zhang \cite{Z}.
\begin{thm}\label{mainkak}
If the subspaces $V_1,\hdots, V_m$ and exponents $\mathbf{p}=(p_1,\hdots,p_m)$ satisfy 
\begin{equation}\label{kscal}
n=\sum_{j=1}^m\frac{\codim(V_j)}{p_j}
\end{equation}
and 
\begin{equation}\label{Vkak}
\dim(V)\leq\sum_{j=1}^m\frac{\dim(V)-\dim(V\cap V_j)}{p_j}
\end{equation}
for all subspaces $V$ of $\mathbb{R}^n$, then
\begin{equation}\label{MKBL}
\int_{\mathbb{R}^n}\prod_{j=1}^m\Biggl(\sum_{T_j\in\mathbb{T}_j}\chi_{T_j}\Biggr)^{1/p_j}\lesssim\delta^n\prod_{j=1}^m(\#\mathbb{T}_j)^{1/p_j}.
\end{equation}
\end{thm}
As we have already alluded to, the inequality \eqref{MKBL} may be viewed as a certain perturbation of  the classical Brascamp--Lieb inequality (this time in the equivalent form \eqref{BLL}). To see this consider the special case where, for each $1\leq j\leq m$, all elements $T_j\in\mathbb{T}_j$ are parallel to $V_j$, so that
$$
\sum_{T_j\in\mathbb{T}_j}\chi_{T_j}=f_j\circ L_j,
$$
where $f_j$ is a sum of characteristic functions of $O(\delta)$ balls in $\mathbb{R}^{n_j}$, and $L_j:\mathbb{R}^n\rightarrow\mathbb{R}^{n_j}$ is a linear surjection with kernel $V_j$. Substituting this into \eqref{MKBL}, and using the scaling condition \eqref{kscal}, we obtain \eqref{BLL} with $r_j=1/p_j$.\footnote{Strictly speaking we obtain \eqref{BLL} only for functions $f_j$ being sums of characteristic functions of $\delta$-balls, although since the implied constant is uniform in $\delta$, we may drop this requirement by a scaling and limiting argument.}  Here we have reconciled the conditions \eqref{Vkak} and \eqref{scaV} via the elementary identity 
$\dim(L_jV)=\dim(V)-\dim(V\cap \ker(L_j)).$ We refer to Maldague \cite{Maldague} for a recent local version of Theorem \ref{mainkak} which has the advantage of not requiring the scaling condition \eqref{kscal}. Kakeya--Brascamp--Lieb inequalities of this type play an important role in the theory of decoupling (also known as Wolff inequalities) -- see \cite{BDAnnals}.

Variants and generalisations of such Kakeya--Brascamp--Lieb inequalities have also been useful in establishing estimates for linear operators. Notably, variants living at intermediate levels of multilinearity, which exploit both curvature and transversality, have led to progress on the original Kakeya maximal conjecture -- see \cite{HRZ} and \cite{Zahl}. Further, certain invariant generalisations of Theorem \ref{mainkak}, where the affine subspaces are also replaced with algebraic varieties (see \cite{Z} and \cite{ZK}) have recently had applications to the $L^p$-improving properties of Radon-like transforms -- see Gressman \cite{Gress}.

\subsection{A Fourier extension variant}\label{Sect:ext}
The Brascamp--Lieb inequality provides a natural framework for certain ``multilinear" variants of the classical Fourier restriction conjecture described in the introduction. In particular, the inequality (E) with $\varepsilon=0$, which is equivalent to the global extension estimate
\begin{equation}\label{npe}
\int_{\mathbb{R}^n}\prod_{j=1}^m|\widehat{g_jd\sigma_j}|^{2/p_j'}\lesssim\prod_{j=1}^m\|g_j\|_{L^2(S_j)}^{2/p_j'},
\end{equation}
may be viewed as a certain oscillatory form of \eqref{BL}. In order to explain this it is convenient to parametrise our (compact) submanifolds of $\mathbb{R}^n$ by smooth mappings $\Sigma_j:U_j\rightarrow\mathbb{R}^n$, where $U_j\subseteq\mathbb{R}^{n_j}$, and define the (parametrised) extension operators
$$
\mathcal{E}_jf_j(x)=\int_{U_j}e^{i\langle x,\Sigma_j(\xi)\rangle}f_j(\xi)d\xi,\;\;\;x\in\mathbb{R}^n.$$
With this notation it is straightforward to verify that \eqref{npe} is equivalent to 
\begin{equation}\label{e}
\int_{\mathbb{R}^n}\prod_{j=1}^m|\mathcal{E}_jf_j|^{2/p_j'}\lesssim\prod_{j=1}^m\|f_j\|_{L^2(U_j)}^{2/p_j'}.
\end{equation}
Now, if $\Sigma_j$ is \emph{linear}, then $\mathcal{E}_jf_j(x)=\widehat{f}_j(\Sigma_j^*(x))$, and so \eqref{e} reduces to
\begin{equation}\label{el}
\int_{\mathbb{R}^n}\prod_{j=1}^m|\widehat{f}_j(\Sigma_j^*(x))|^{2/p_j'}dx\lesssim\prod_{j=1}^m\|f_j\|_{L^2(U_j)}^{2/p_j'},
\end{equation}
which by Plancherel's theorem applied in $\mathbb{R}^{n_j}$ for each $j$, becomes
\begin{equation}\label{elp}
\int_{\mathbb{R}^n}\prod_{j=1}^m g_j(\Sigma_j^*(x))^{1/p_j'}dx\lesssim\prod_{j=1}^m\Bigl(\int_{\mathbb{R}^{n_j}}g_j\Bigr)^{1/p_j'}.
\end{equation} 
This is of course the Brascamp--Lieb inequality in the equivalent form \eqref{BLL} with $r_j=1/p_j'$ and $L_j=\Sigma_j^*$. In this way, permitting nonlinear $\Sigma_j$ in \eqref{e}, or equivalently, nonlinear submanifolds $S_j$ in \eqref{npe}, constitutes a generalisation of the classical Brascamp--Lieb inequality. It should be noticed that unless $p_1=\cdots=p_m$, such inequalities are not (manifestly at least) bounds on multilinear operators, and differ in that respect from the linear and nonlinear Brascamp--Lieb inequalities discussed earlier.

At this level of generality we have the following result from \cite{BBFL}. For convenience we suppose that $0\in U_j$ for each $1\leq j\leq m$, which is of course without loss of generality.
\begin{thm}[\cite{BBFL}] \label{mainext}
If $\BL(d\mathbf{\Sigma}(0)^*,\mathbf{p}')<\infty$ then provided the $U_j$ are sufficiently small neighbourhoods of $0$,
\begin{equation}\label{ef}
\int_{B_R}\prod_{j=1}^m|\mathcal{E}_jf_j|^{2/p_j'}\lesssim_\varepsilon R^\varepsilon\prod_{j=1}^m\|f_j\|_{L^2(U_j)}^{2/p_j'}.
\end{equation}
Here $\mbox{d}\mathbf{\Sigma}(0)^*$ denotes the $m$-tuple of linear maps $(d\Sigma_1(0)^*,\hdots,d\Sigma_m(0)^*)$.
\end{thm}
As mentioned in the introduction, this result with Loomis--Whitney data, that is, where the maps $\mathbf{L}=\mbox{d}\mathbf{\Sigma}(0)^*$ satisfy the basis condition \eqref{LoomW}, is referred to as the multilinear restriction inequality, and originates in \cite{BCT}. 
Estimates of this type have proved to have many applications, most notably to the classical (linear) restriction conjecture -- see \cite{BG} -- and to the theory of decouplings (also known as Wolff inequalities) -- see \cite{BDAnnals}.
%Section \ref{Sect:App} for further discussion.
Remarkably, Theorem \ref{mainext} and (the Kakeya--Brascamp--Lieb) Theorem \ref{mainkak} are virtually equivalent -- see \cite{BBFL}; we refer to Bejenaru \cite{Bej} for an alternative approach to multilinear extension estimates which avoids explicit Kakeya-type considerations.

\subsection{A multilinear oscillatory integral variant}\label{Sect:mosc}
Another, rather different oscillatory form of the Brascamp--Lieb inequality was introduced by Christ--Li--Tao--Thiele \cite{CLTT}, this time exhibiting connections with questions related to Szemer\'edi's theorem from additive combinatorics. Motivated by questions of boundedness of certain multilinear oscillatory singular integral operators, the authors consider the multilinear functional
\begin{equation}\label{clttfun}
\Lambda_\lambda(f_1,\hdots,f_m)=\int_{\mathbb{R}^n}e^{i\lambda P(x)}\prod_{j=1}^m f_j(L_jx)\eta(x)\,dx,
\end{equation}
where $P$ is a real-valued measurable function, $\eta$ is a smooth bump function, and $\lambda$ is a real parameter. Of course one has the elementary estimate
\begin{equation*}
|\Lambda_\lambda(f_1,\hdots,f_m)|\lesssim\prod_{j=1}^m\|f_j\|_{\infty}
\end{equation*}
regardless of the linear maps $L_j$ and phase function $P$. In this context it is of particular interest to identify conditions on the data $(\mathbf{L}, P)$ for which there is some additional decay in $\lambda$ -- that is, for which the right hand side above may include an additional factor of $\lambda^{-\varepsilon}$ for some 
$\varepsilon>0$. Such $(\mathbf{L},P)$  are said to have the \emph{power decay property}, and the objective in \cite{CLTT}, and subsequent works, is to characterise these.
Of course, this set-up naturally generalises to incorporate $L^{p_j}(\mathbb{R}^{n_j})$ norms of $f_j$ for each $j$, allowing the theory to generalise (or interact with) that of the classical Brascamp--Lieb inequality \eqref{BL} -- see \cite{BCCT2}. 

The function $P$ is said to be degenerate relative to $\mathbf{L}$ if it may be expressed as a linear combination of measurable functions of the form $\phi_j\circ L_j$. It is straightforward to see that such data cannot have the power decay property, and the objective is to establish power decay in all other situations. Although this remains a challenging open problem, we refer the reader to \cite{Chr, ChristDiogo, GGX, GOX, GU, GX, NOZ, X} for substantial recent progress in this direction.

\subsection{A singular integral variant}
At the level of examples at least, singular integral variants of the Brascamp--Lieb inequality \eqref{BL} have been the focus of considerable attention over the last two decades, including the celebrated boundedness of the bilinear Hilbert transform of Lacey and Thiele \cite{LT} and subsequent developments (see, for example, Demeter--Pramanik--Thiele \cite{DPT} and Muscalu--Tao--Thiele \cite{MTT}). In the recent survey article of Durcik--Thiele \cite{DT} a general set-up is presented that differs from \eqref{BL} in that some of the functions $f_j$ corresponding to exponents $p_j=1$ are taken to be Calder\'on--Zygmund singular integral kernels, rather than Lebesgue integrable functions. Thus inequalities of the form
\begin{equation}\label{SBL}
\Biggl|\int_{\mathbb{R}^n}\prod_{j=1}^kf_j(L_jx)\prod_{j'=k+1}^mK_j(L_jx)\,dx\Biggr|\lesssim\prod_{j=1}^k\|f_j\|_{L^{p_j}(\mathbb{R}^{n_j})}
\end{equation}
are sought, where for $k+1\leq j\leq m$, the distribution $K_j$ is a standard Calder\'on--Zygmund kernel. A well-known example is the conjectural bound on the so-called triangular Hilbert transform, given by
$$
\Biggl|\int_{\mathbb{R}^3} f_1(x_2,x_3)f_2(x_1,x_3)f_3(x_1,x_2)\frac{dx}{x_1+x_2+x_3}\Biggr|\lesssim\|f_1\|_{L^3(\mathbb{R}^2)}\|f_2\|_{L^3(\mathbb{R}^2)}\|f_3\|_{L^3(\mathbb{R}^2)},
$$
which amounts to a singular Brascamp--Lieb inequality where the linear maps are the three-dimensional Loomis--Whitney maps with a fourth $L_4x:=x_1+x_2+x_3$ appended. Here the singular kernel is the classical Hilbert principal value kernel, and the integral is interpreted accordingly. As is pointed out in \cite{DT}, certain aspects of the Brascamp--Lieb theory already discussed usefully apply in this singular context, such as the Fourier-invariance property \eqref{inv}. Furthermore, since delta distributions are examples of singular kernels, the finiteness conditions for \eqref{BL} (\eqref{scapf} and \eqref{scapfV} with $p_{k+1}=\cdots=p_m=1$) also become necessary conditions for finiteness in \eqref{SBL}. As may be expected, the theory of \eqref{SBL} is much less developed and largely very different from that of \eqref{BL}, involving time-frequency analysis and other methods that exploit cancellation. We refer to \cite{DT, DT2, MZ} for further discussion and recent results.

\subsection{A variant on LCA groups}
The Brascamp--Lieb inequality is also naturally formulated in the context of locally compact abelian groups, whereby \eqref{BL} becomes
\begin{equation}\label{BLLCA}
\int_G \prod_{j=1}^m f_j\circ\varphi_j\leq C\prod_{j=1}^m\|f_j\|_{L^{p_j}(G_j)},
\end{equation}
where the maps $\varphi_j:G\rightarrow G_j$ are homomorphisms of LCA groups, and the integrals involved are with respect to (suitably normalised) Haar measures. This departs from the euclidean theory in some interesting ways. In the context of finitely generated discrete groups a finiteness characterisation was established in \cite{BCCT2}. Motivated by problems in communication theory, this was studied further by 
Christ--Demmel--Knight--Scanlon--Yelick \cite{CDKSY}, who showed that a polynomial time algorithm for verifying the finiteness condition is equivalent to Hilbert's tenth problem over the rational numbers. Furthermore, for torsion-free discrete groups, they showed that if the constant $C$ is finite then it must equal $1$ -- a conclusion consistent with the classical Young's convolution inequality on the integers. Similar results may be established in the compact setting, and ultimately give rise to an abstract duality principle of the form \eqref{inv} -- see \cite{BJ} for  more general conclusions and clarification.

\subsection{A measure space variant}
Continuing somewhat in the spirit of the preceding formulation, it is reasonable to consider variants of the Brascamp--Lieb inequality on more general ambient spaces. Whilst the classical form of the Brascamp--Lieb inequality on euclidean spaces is amenable to a variety of techniques, the approach based on a diffusion flow has proved to be particularly effective in more general contexts. An elegant example in the setting of real spheres takes the form
\begin{equation} \label{CLL_sphere}
\int_{\mathbb{S}^{n-1}} \prod_{j=1}^n f_j(x_j) d\sigma(x) \leq \prod_{j=1}^n \Bigl(\int_{\mathbb{S}^{n-1}} f_j(x_j)^2 d\sigma(x) \Bigr)^{1/2}.
\end{equation} 
This sharp inequality was derived in \cite{CLL} by using a heat flow monotonicity argument for functions on $\mathbb{S}^{n-1}$. In this case, sharpness refers both to the optimality of the constant $1$ and also in the sense of the $L^2$ norms on the right hand side (which provides a clear advantage over a simple application of the H\"older inequality). Despite the fact that the underlying Brascamp--Lieb mappings are simply coordinate projections, obviously the coordinates of points on spheres are not independent and thus there is no Fubini-type identity on $L^1$. The inequality \eqref{CLL_sphere} may be viewed as a correlation inequality which quantifies how far coordinate functions on spheres are from being independent.

Carlen, Lieb and Loss applied similar heat flow arguments in the context of the classical Brascamp--Lieb inequality \eqref{BL} for rank-one mappings (also in \cite{CLL}) and in the seemingly rather different setting of permutation groups in order to obtain a version of Hadamard's inequality for matrix permanents (see \cite{CLL2}). A unification and significant extension may be found in \cite{BCLM} (see also the precursor \cite{BCM}) whose setting consisted of an underlying measure space, a Markov semigroup with generator $L$ acting on the measure space, and Brascamp--Lieb mappings between measure spaces which interact in an appropriate manner with the generator $L$. Brascamp--Lieb type inequalities in such a setting were generated through an abstract argument using the semigroup interpolation method and this facilitated a comprehensive understanding of several concrete settings. Further more recent developments along these lines may be found in work of Bramati \cite{Bramati}.

\subsection{A Lorentz norm variant}
Whilst the variants discussed up to this point have amounted to modification to the form on the left hand side of the classical Brascamp--Lieb inequality, it is also beneficial to extend the framework to allow input functions in Lorentz spaces rather than the standard Lebesgue spaces. One need look no further than the classical and ubiquitous Hardy--Littlewood--Sobolev inequality (in dual form) for a special case. More recently though more elaborate Brascamp--Lieb type inequalities such as
\begin{equation*}
\int_{\mathbb{C}^{2k+1}} \frac{f_0(x_0 - x_1 + \cdots - x_{2k-1} + x_{2k}) \prod_{j=1}^{2k} f_j(x_j)}{|x_0 - x_1||x_1 - x_2| \cdots |x_{2k-1} - x_{2k}|} \, d\mu(x_0,\ldots,x_{2k}) \lesssim  \prod_{j=0}^{2k} \|f_j\|_2, 
\end{equation*}
where $d\mu$ is product measure, have arisen in inverse scattering theory -- see work of Brown \cite{Brown}, as well as \cite{AFR, BOP,Perry}. As observed by Christ in \cite{Perry}, thanks to multilinear interpolation it is possible to obtain estimates of the form\footnote{Here we are using the standard notation $L^{p,r}$ for Lorentz spaces.}
\begin{equation} \label{BL_Lorentz}
\int_{\R^n}\prod_{j=1}^m f_j(L_jx) \,dx
\lesssim \prod_{j=1}^m\|f_j\|_{L^{p_j,r_j}(\R^{n_j})}
\end{equation}
from the finiteness of $\BL(\datum{L},\datum{\widetilde{p}})$, for $\datum{\widetilde{p}}$ in a neighbourhood of $\datum{p}$, and under the constraint $\sum_{j=1}^m 1/r_j \geq 1$. This generates, for instance, the estimate \eqref{BL_Lorentz} in the case of simple\footnote{A Brascamp--Lieb data $(\datum{L},\datum{p})$ is said to be \emph{simple} if the scaling condition \eqref{sca} holds and the dimension condition \eqref{scaV} holds with strict inequality for all proper and non-trivial subspaces.} Brascamp--Lieb data $(\datum{L},\datum{p})$ and whenever $\sum_{j=1}^m 1/r_j \geq 1$. This latter restriction on the exponents $\datum{r}$ has been shown to be necessary in the case of arbitrary Brascamp--Lieb data (see Bez--Lee--Nakamura--Sawano \cite{BLNS}). A characterisation of allowable estimates of the form \eqref{BL_Lorentz} in the case of non-simple data appears to be an interesting open problem\footnote{It was observed in \cite{BLNS} that the stronger condition $\sum_{j=1}^n 1/r_j \geq \frac{n}{n-1}$ is necessary for Loomis--Whitney data (with one-dimensional kernels), so we may expect such a characterisation to be somewhat subtle.}.

The ``upgrading" trick of Christ has also been employed recently in work of Kato--Miyachi--Tomita \cite{KMT} on the boundedness of certain bilinear pseudodifferential operators. Naturally appearing in their study are bilinear symbol classes associated to weights $\widetilde{V}$ defined on $\mathbb{R}^n \times \mathbb{R}^n$ for which the weighted (discrete) Brascamp--Lieb inequality 
\[
\sum_{\ell_1,\ell_2 \in \mathbb{Z}^n} A_1(\ell_1 + \ell_2)A_2(\ell_1)A_3(\ell_2)V(\ell_1,\ell_2) \lesssim \|A_1\|_{\ell^2(\mathbb{Z}^n)} \|A_2\|_{\ell^2(\mathbb{Z}^n)} \|A_3\|_{\ell^2(\mathbb{Z}^n)}
\]
holds. Here, $V$ is a discretisation of $\widetilde{V}$, and the observation of Christ was used to show that $V \in \ell^{4,\infty}(\mathbb{Z}^n \times \mathbb{Z}^n)$ are admissible.

\subsection{Some remarks on connections between the different variants}
While quite different in many respects, the variant Brascamp--Lieb inequalities discussed above have interrelations on a number of different levels.  As may be expected from classical (linear) restriction theory, when the $S_j$ are hypersurfaces the extension variant is easily seen to imply the Kakeya-type variant via a routine randomisation argument involving wavepackets -- see, for example, \cite{BCT} or \cite{BBFL}. Furthermore, as mentioned at the end of Section \ref{Sect:kak}, the extension variant may be deduced from the Kakeya variant by a routine induction-on-scales argument, provided one is prepared to accept a small power or polylogarithmic loss in the truncation parameter $R$. A link between the nonlinear and extension variants becomes apparent on extending \eqref{e} to incorporate more general H\"ormander-type oscillatory integral operators -- an exercise that leads to a common generalisation of the two (see \cite{B}). The nonlinear, extension and Kakeya-type variants are also related from a methodological point of view, all being accessible (up to endpoints at least) by the method of induction-on-scales -- see \cite{BBFL} for further details.

A  somewhat superficial connection between the nonlinear and the oscillatory variants becomes apparent upon expressing one of the functions (say $f_m$) in the nonlinear setting in terms of its Fourier transform using the inversion formula, upon which the map $B_m$ essentially emerges as the phase $P$. A less superficial connection between the extension, singular integral and nonlinear variants stems from the simple observation that the three-dimensional (conjectured) endpoint multilinear extension inequality of \cite{BCT}, namely
$$
\|\widehat{g_1d\sigma_1}\widehat{g_2d\sigma_2}\widehat{g_3d\sigma_3}\|_{L^1(\mathbb{R}^3)}\lesssim\|g_1\|_{L^2(S_1)}\|g_2\|_{L^2(S_2)}\|g_3\|_{L^2(S_3)},
$$
implies
$$
\Bigl|\int_{\mathbb{R}^3}\widehat{g_1d\sigma_1}\widehat{g_2d\sigma_2}\widehat{g_3d\sigma_3}\,\widehat{K}\Bigr|\lesssim\|g_1\|_{L^2(S_1)}\|g_2\|_{L^2(S_2)}\|g_3\|_{L^2(S_3)},$$
for any three-dimensional Calder\'on--Zygmund kernel $K$. Applying Parseval's identity, this is equivalent to
$$
\Bigl|\int_{\mathbb{R}^3}g_1d\sigma_1*g_2d\sigma_2*g_3d\sigma_3 \,K\Bigr|\lesssim\|g_1\|_{L^2(S_1)}\|g_2\|_{L^2(S_2)}\|g_3\|_{L^2(S_3)}$$
holding for any such $K$. If the surfaces $S_1, S_2, S_3$ are linear then this is an (albeit rather straightforward) instance of a singular Brascamp--Lieb inequality. However, for nonlinear $S_j$ it is not known, and appears to be an interesting challenge to existing methods. We note that in the special case where $K$ is the delta distribution at a point, this is an instance of the convolution inequality (C) originating in \cite{BHT}, which is in turn a manifestation of the nonlinear Brascamp--Lieb inequality for Loomis--Whitney data \cite{BCW}.

\section{On the proof of Theorem \ref{tfae}}\label{SectProof}  In this section we show how Theorem \ref{tfae} may be reduced to the wider theory of the Brascamp--Lieb inequality described in Section \ref{SectWider}. To a great extent this is achieved by amalgamating the results of \cite{BCCT2, BBG, BBFL, BBBCF}, with the exception of the implication from (E) to (T), which has some subtlety due to the $\varepsilon$-loss in (E). We recall that the scaling condition \eqref{scaleawayscaleawayscaleaway} should be assumed here throughout.

\subsection{The equivalence of (T) and (E)}
That (T)$\implies$(E) is quickly reduced to Theorem \ref{mainext}, which as we will see, is simply a version of this statement in terms of the parametrised extension operators of Section \ref{Sect:ext}. First of all, since the submanifolds $S_j$ are compact, it is enough to establish a local version of (T)$\implies$(E), whereby the the tangent spaces in (T) are at a single point, and the conclusion (E) is for functions $g_j$ supported in a sufficiently small neighbourhood of that point. For convenience we take this distinguished point to be the origin, assuming as we may that $0\in S_j$ for each $j$, and denote by (T$_0$) and (E$_0$) the conditions (T) and (E) with these restrictions.  Since (E$_0$) and its parametrised version \eqref{ef} are equivalent, the equivalence of (T) and (E) will follow if we show that
\begin{itemize}
\item[(i)]
(T$_0$) is equivalent to the finiteness of $\BL(d\mathbf{\Sigma}(0)^*,\mathbf{p}')$, and 
\item[(ii)] \eqref{ef} for all $\varepsilon>0$ implies the finiteness of $\BL(d\mathbf{\Sigma}(0)^*,\mathbf{p}')$.
\end{itemize}
In doing this, for convenience we shall assume, as we may, that the $S_j$ are parametrised in such a way that $0\in U_j$.

Part (i) follows directly from Theorem \ref{t:fin} along with the elementary observation that
\begin{equation}\label{elobs}
\dim(d\Sigma_j(0)^*V)=\dim(V)-\dim(V\cap (T_0S_j)^\perp)
\end{equation}
for all subspaces $V$ of $\mathbb{R}^n$.

Part (ii) would follow from a scaling and limiting argument were we permitted to take $\varepsilon=0$ in (E). Since we are not, we appeal to an adaptation of the Knapp-type examples from the classical linear and bilinear restriction theory (see \cite{TVV}). To this end let $V$ be a subspace of $\mathbb{R}^n$, and suppose (E$_0$) holds for all $\varepsilon>0$. By Theorem \ref{t:fin} it suffices to show that 
\begin{equation}\label{sufftoshow}
\dim(V)\leq \sum_{j=1}^m\frac{\dim(L_jV)}{p_j'},
\end{equation}
where $L_j:=d\Sigma_j(0)^*$. For each $1\leq j\leq m$ we define a set $X_j\subseteq U_j$ by
$$
X_j= \mathcal{N}_\delta\left((L_jV)^\perp\cap B(0,\delta^{1/2})\right),
$$
where $\mathcal{N}_\delta(E)$ is used to denote the $\delta$-neighbourhood of a set $E$. We clarify that $(L_jV)^\perp$ is a subspace of $\mathbb{R}^{n_j}$, so that $B(0,\delta^{1/2})$ should be interpreted as a $\delta^{1/2}$-ball in $\mathbb{R}^{n_j}$.
Similarly we define a set $X\subseteq\mathbb{R}^n$ by 
$$
X=\mathcal{N}_{c\delta^{-1/2}}\left(V\cap B(0,c\delta^{-1})\right)
$$
for some suitably small constant $c>0$, independent of $\delta$. Setting $f_j=\chi_{X_j}$ we have
$$
\mathcal{E}_jf_j(x)=\int_{X_j}e^{ih_j(x,\xi)}d\xi,
$$
where 
\begin{equation}\label{plumberjustbeen}
h_j(x,\xi):=\langle x,\Sigma_j(\xi)\rangle=\langle L_jx,\xi\rangle+\langle x,\Sigma_j(\xi)-d\Sigma_j(0)\xi\rangle.
\end{equation}
Now, we claim that $|h_j(x,\xi)|\leq 1/5$ whenever $\xi\in X_j$ and $x\in X$, provided the constant $c$ is taken to be sufficiently small. 
Using that fact we have that $|\mathcal{E}_jf_j(x)|\gtrsim |X_j|$ for all $x\in X$. Hence by (E$_0$) with $R=\delta^{-1}$, 
$$
|X|\prod_{j=1}^m|X_j|^{1/p_j'}\lesssim_\varepsilon \delta^{-\varepsilon}
$$
for all $0<\delta\ll 1$ and all $\varepsilon>0$. Since
$$
|X_j| \sim (\delta^{1/2})^{n_j-\dim(L_jV)}\delta^{\dim(L_jV)}
$$
and
$$
|X| \sim (\delta^{-1})^{\dim(V)}(\delta^{-1/2})^{n-\dim(V)},
$$
\eqref{sufftoshow} follows by taking both $\varepsilon$ and $\delta$ arbitrarily small.

It remains to prove the claim, which by \eqref{plumberjustbeen}, will follow if we can show that 
\begin{equation}\label{split}
|\langle L_jx,\xi\rangle|,\,|\langle x,\Sigma_j(\xi)-d\Sigma_j(0)\xi\rangle|\leq 1/10
\end{equation}
whenever $\xi\in X_j$ and $x\in X$, again provided the constant $c$ is taken to be sufficiently small. However, these are elementary exercises using the basic geometries of the sets $X_j$ and $X$ for the first, and the $C^2$ character of the maps $\Sigma_j$ for the second. We leave these to the interested reader.

\subsection{The equivalence of (T) and (C)} The arguments presented here are largely re-workings of arguments from \cite{BBBCF}.
Merely for convenience we begin by expressing (C) in a parametrised form using the smooth maps $\Sigma_j:U_j\rightarrow\mathbb{R}^n$ defined in Section \ref{Sect:ext}, so that surface measure $d\sigma_j$ on $S_j$ is given
by 
$$
\int_{S_j}fd\sigma_j=\int_{U_j}f(\Sigma_j(x))J_j(x)dx.
$$
Here, by compactness, the jacobian factor $J_j(x)$ is both bounded from above and below on the parameter set $U_j\subseteq\mathbb{R}^{n_j}$, where $n_j=\dim(S_j)$. Next we observe that for any $a\in\mathbb{R}^n$,
\begin{eqnarray}\label{e:multconvatx}
\begin{aligned}
g_1d\maybed\sigma_1 * \dots * g_md\maybed\sigma_m(a)
&= \int_{U_1 \times \dots \times U_m}  f_1 \otimes \dots \otimes f_m(y) \delta(F_a(y))  J(y)\,dy\\
&\sim \int_{U_1 \times \dots \times U_m}  f_1 \otimes \dots \otimes f_m(y) \delta(F_a(y))  \,dy\end{aligned}
\end{eqnarray}
where $J=J_1\otimes\cdots\otimes J_m$, $f_j = g_j \circ \Sigma_j$, and
\begin{equation} \label{e:convF}
F_a(y) = \Sigma_1(y_1) + \dots + \Sigma_m(y_m) - a.
\end{equation}
Consequently, (C) reduces to the nonlinear Brascamp--Lieb inequality
\begin{equation}\label{redu}
\int_{M} f_1\otimes\cdots\otimes f_m d\mu_M\lesssim\prod_{j=1}^m\|f_j\|_{L^{p_j}(U_j)},
\end{equation}
where the manifold $$M= \{ (y_1,\ldots,y_m) \in U_1 \times \dots \times U_m : F_a(y) = 0\}. $$
Of course in order to have the required $L^\infty$ bound on the iterated convolution in \eqref{e:multconvatx}, the bound in \eqref{redu} should be  seen to be locally uniform in the implicit parameter $a$. 
In this case the compactness of the surfaces $S_j$ allows us to sidestep this point, and without loss of generality (using the translation-invariance of our hypotheses) suppose that $a=0$, and write $F_0=F$. Similarly for convenience we assume, as we may,  that $\Sigma_j(0)=0$, so that in particular $0$ is a point on $M$.

Of course the main implication to establish is (T)$\implies$(C), for which we appeal to the nonlinear Brascamp--Lieb inequality of Theorem \ref{main'}.
By Theorem \ref{main'}  it suffices to prove that $\BL(T_0M,\datum{p}) < \infty$, where
$
T_0M
$
is the tangent space to $M$ at the origin. For this it remains to see that the Brascamp--Lieb finiteness condition  \eqref{scapfV} for the data $(T_0M, \mathbf{p})$ is equivalent to (T). While this may be established directly with some careful linear algebra, it is more convenient to first appeal to the Fourier duality relation \eqref{inv} to replace the data $(T_0M, \mathbf{p})$ here with its dual data $((T_0M)^\perp,\datum{p}')$. Since
\[
(T_0M)^\perp = \{ (d\Sigma_1(0)^*x,\ldots,d\Sigma_m(0)^*x) : x \in \R^n \},
\]
the finiteness condition \eqref{scapfV} for the data $((T_0M)^\perp,\datum{p}')$, and hence $(T_0M, \mathbf{p})$, follows quickly from (T) using the elementary identity \eqref{elobs}.

Finally, the converse (C)$\implies$(T) is a direct consequence of the following minor variant of Lemma 6 of \cite{BBG}.
\begin{prp}
Suppose that $M$ is a smooth submanifold of  $\mathbb{R}^{n_1}\times\cdots\times\mathbb{R}^{n_m}$, and $\mathbf{p}=(p_1,\hdots,p_m)\in [1,\infty)^m$ is such that
$$\dim(M)=\sum_{j=1}^m\frac{n_j}{p_j}$$ and
$$
\int_M f_1\otimes\cdots\otimes f_m\,d\mu_M\leq C\prod_{j=1}^m\|f_j\|_{L^{p_j}(\mathbb{R}^{n_j})}
$$
for some constant $C$, where $d\mu_M$ denotes surface measure on $M$. Then $\BL(T_xM,\mathbf{p})\leq C$ for all $x\in M$.
\end{prp}

\section{Applications}\label{Sect:App}
The estimates (C) and (E) in the statement of Theorem \ref{tfae} have found numerous applications. Most (although not all -- see \cite{BBBCF} for applications in abstract harmonic analysis) stem from the close relationship between transversality and curvature, and the well-known fact that harmonic analysis is a powerful tool for understanding many analytical problems where some underlying curved manifold plays an important role. In this section we elaborate on certain (PDE) applications of the lesser-known estimate (C), referring the reader to \cite{B}, \cite{Demeter} for further discussion of (E).

In the analysis of nonlinear dispersive equations, estimates involving iterated convolutions frequently arise.
A widely used technique (see, for example, Beals \cite{Beals}, Bourgain \cite{Bourgain93},  Klainerman--Machedon \cite{KM}) in the study of the local well-posedness of subcritical problems is the use of $X^{s,b}_S$ spaces, associated with an underlying surface $S$, in an iteration argument.
In this way, it is desirable to establish control on the nonlinearity in the setting of these function spaces and, under certain structural assumptions on the nonlinearity, it is often possible to reduce estimates of this type to multilinear convolution estimates.
Naturally it is desirable to have a general framework for such estimates, and a systematic study of weighted convolution estimates of the form
\begin{equation} \label{e:Taoweighted}
\bigg| \int_{G} f_1 \otimes \dots \otimes f_m(y) \, w(y) \,d\mu_G(y) \bigg|
\leq C \prod_{j=1}^m \|f_j\|_{L^2(Z)}
\end{equation}
for $L^2$ functions was carried out by Tao \cite{TaoAJM}.
Here, $Z$ is an abelian group, written additively, $G$ is the subgroup $\{ y \in Z^m : y_1 + \dots + y_m = 0\}$ of $Z^m$, and $\mu_G$ is a Haar measure and $w$ is a (weight) function on $G$.
We refer the reader to \cite{TaoAJM} for further details concerning the motivation for estimates of the type \eqref{e:Taoweighted}, along with applications to bilinear estimates associated with the KdV, wave and Schr\"odinger equations.

The multilinear convolution estimate (C), which involves singular measures (rather than weights), falls outside the scope of \eqref{e:Taoweighted}. However, such singular convolution estimates also arise naturally in dispersive PDE, and in particular, have played a role in recent breakthroughs in the low-regularity theory of the Zakharov system by Bejenaru--Herr--Holmer--Tataru \cite{BHHT} and Bejenaru--Herr \cite{BH}.  This is a coupled system of nonlinear Schr\"odinger and wave equations, and takes the form
\begin{align*}
(i\partial_t + \Delta)u & = nu \\
(\partial_t^2 - \Delta)n & = \Delta |u|^2,
\end{align*}
for functions $u : \R^{d+1} \to \mathbb{C}$ and $n : \R^{d+1} \to \R$. Formulated by Zakharov \cite{Zakharov}, this arises as a model in plasma physics, and much effort has been spent developing the well-posedness theory of this system (see, for example, \cite{BH,BHHT,BC,GTV,OT,SS}). For initial data in $L^2$-based Sobolev spaces, the papers \cite{BHHT} and \cite{BH} provided definitive results on the local well-posedness of the system in two and three spatial dimensions, respectively, in the subcritical regime.
The iteration argument in \cite{BHHT} involves multilinear estimates in terms of certain $X^{s,b,1}_{S}$ spaces, where $S$ is either a paraboloid or a cone in $\R^3$.
The underlying qualitative phenomenon is the fact that, given three $C^2$, bounded, and transverse\footnote{In this context, if $\nu_1, \nu_2, \nu_3$ are linearly independent for any choice of normal vector $\nu_j$ to $S_j$, we say that the $S_j$ are \emph{transverse}.} surfaces $S_1$, $S_2$ and $S_3$ in  $\R^3$, the convolution of two $L^2$ functions supported on two of the surfaces has a well-defined restriction, as an $L^2$ function, to the third surface $S_3$.
This phenomenon, which has already been touched on in the introduction, is captured by the estimate
\begin{equation*} \label{e:BHT}
\|g_1d\maybed\sigma_1 * g_2d\maybed\sigma_2\|_{L^2(\maybedsigma_3)} \leq C \|g_1\|_{L^2(\maybedsigma_1)} \|g_2\|_{L^2(\maybedsigma_2)},
\end{equation*}
or equivalently, by duality,
\begin{equation} \label{e:BHTdual}
|g_1d\maybed\sigma_1 * g_2d\maybed\sigma_2 * g_3d\maybed\sigma_3(0)| \leq C \prod_{j=1}^3 \|g_j\|_{L^2(\maybedsigma_j)}.
\end{equation}
This is of course a particular instance of (C), and the associated 
transversality assumption (T) reduces to the three-dimensional case of \eqref{lwt}.
The specific bound \eqref{e:BHTdual}, in the quantitative form of relevance to applications, was first established in \cite{BHT} (see also \cite{BCW}).
We remark that in most applications of this type, the underlying submanifolds $S_j$ satisfy such a ``Loomis--Whitney-type" transversality condition, whereby (T) reduces to the ``basis condition"
$$
(TS_1)^{\perp}\oplus\cdots\oplus (TS_m)^\perp=\mathbb{R}^n,
$$
where (necessarily) $p_1=\cdots=p_m=m-1$; cf. \eqref{LoomW}. The reader will of course recognise this as a generalisation of \eqref{lwt} from the introduction, and we refer to \cite{BB} and \cite{BBBCF} for further discussion of its significance.

There are further, quite varied applications of (C), and the nonlinear Brascamp--Lieb inequality more generally, in PDE-related problems.
For example, in their study of the weakly nonlinear large-box limit of the two-dimensional cubic nonlinear Schr\"odinger equation, Faou--Germain--Hani \cite{FGH} derive a new nonlinear integro-differential equation governed by a trilinear form that draws on the nonlinear Loomis--Whitney inequality from \cite{BCW}.
Also, in the spirit of \cite{BHHT} and \cite{BH}, nonlinear Loomis--Whitney-type inequalities have been used very recently to push forward the mathematical theory of the Klein--Gordon--Zakharov system in two dimensions \cite{K} and a system of quadratic derivative nonlinear Schr\"odinger equations \cite{HK}.
In somewhat different territory, explicit examples of multilinear Radon-like transforms of the form \eqref{functionalanalysis} have also appeared in obstacle scattering, in particular, in the recovery of singularities of a potential $q$ by its so-called Born series $q_B$; we refer the reader to \cite{RV} (see also \cite{BBG2}) for further background and discussion on the manner in which bounds of the form \eqref{great!} naturally arise in this context.
It seems reasonable to expect further such applications arising from Theorems \ref{tfae} and \ref{main'}.

\begin{ack}
We would like to thank Satoshi Masaki and Hideo Takaoka for their generous hospitality at RIMS in the summer of 2019, where some of this work was carried out. We would also like to thank Stefan Buschenhenke, Michael Cowling and Taryn Flock for helpful discussions during the preparation of this article. Finally, we thank the anonymous referee for their helpful suggestions.
\end{ack}

\end{document}